\input amssym.def
\input amssym
\magnification=1200
\parindent0pt
\hsize=16 true cm \baselineskip=13  pt plus .2pt $ $

\def\Z{\Bbb Z}
\def\A{\Bbb A}
\def\S{\Bbb S}
\def\D{\Bbb D}
\def\R{\Bbb R}
\def\P{\Bbb P}

\centerline {\bf  On finite groups acting on spheres}

\centerline {\bf  and finite subgroups of orthogonal groups}

\bigskip \bigskip

\centerline {Bruno P. Zimmermann}

\bigskip

\centerline {Universit\`a degli Studi di Trieste}
\centerline {Dipartimento di Matematica e Informatica}
\centerline {34127 Trieste, Italy}

\bigskip \bigskip

Abstract.  {\sl  This is a survey on old and new results as well as an
introduction to various related basic notions and concepts, based on two talks
given at the International Workshop on Geometry and Analysis in Kemerovo
(Sobolev Institute of Mathematics, Kemerovo State University) and at the
University of Krasnojarsk in June 2011. We discuss finite groups acting on
low-dimensional spheres, comparing with the finite subgroups of the
corresponding orthogonal groups, and also finite simple groups acting on
spheres and homology spheres of arbitrary dimension.}

\bigskip \bigskip

{\bf 1. Introduction}

\medskip

We are interested in the class of finite groups which admit an
orientation-preserving action on a sphere $S^n$ of a given dimension $n$.  All
actions in the present paper will be faithful and orientation-preserving (but
no necessarily free). So formally an action of a finite group $G$ is an
injective homomorphism from $G$ into the group of orientation-preserving
homeomorphism of $S^n$; informally, we will consider $G$ as a group of
orientation-preserving homeomorphisms of $S^n$ and distinguish various types of
actions:

\bigskip

{\it topological actions}: $G$ acts by homeomorphisms;

\medskip

{\it smooth actions}: $G$ acts by diffeomorphisms;

\medskip

{\it linear (orthogonal) actions}: $G$ acts by orthogonal maps of $S^n \subset \Bbb
R^{n+1}$ (that is, $G$ is a subgroup of the orthogonal group ${\rm SO}(n+1)$); or, more
generally, any topological action which is conjugate to a linear action;

\medskip

{\it locally linear actions}: topological actions which are linear in regular
neighbourhoods of fixed points of any nontrivial subgroup.

\medskip

It is well-known that smooth actions are locally linear (by the existence of
equivariant regular neighbourhoods, see [Bre]; see also the discussion in the next
two  sections for some examples).

\medskip

The reference model for finite group actions on $S^n$ is the orientation-preserving
orthogonal group ${\rm SO}(n+1)$; in fact, the only examples of actions which can
easily be seen are linear actions by finite subgroups of ${\rm SO}(n+1)$. A rough
general guiding line is then the following:

\bigskip

{\it Motivating naive conjecture}: Every action of a finite group on $S^n$ is
linear (that is, conjugate to a linear action).  More generally, how far can an
action by homeomorphisms or diffeomorphisms be from a linear action?

\medskip

It is one of the main features of linear actions that fixed point sets of
single elements  are standard unknotted spheres $S^k$ in $S^n$ (that is, intersections
of linear subspaces $\Bbb R^{k+1}$ of $\Bbb R^{n+1}$ with $S^n$). We note that it is a
classical and central result of Smith fixed point theory that, for an action
of a finite $p$-group (a group whose order is a power of a prime $p$) on a mod
$p$ homology $n$-sphere (a closed $n$-manifold with the mod $p$ homology of
$S^n$), fixed point sets are again mod $p$ homology spheres (see [Bre]; one has to
consider homology with coefficients in the integers mod
$p$ here since this does not remain true in the setting of integer homology
spheres).

\medskip

We note that for {\it free actions} on $S^n$ (nontrivial elements have empty
fixed point sets), the class of finite groups occurring is very restricted
(they have periodic cohomology, of period $n+1$, see [Bro]). On the other hand,
for not necessarily free actions, all finite groups occur for some $n$  (by
just considering a  faithful, real, linear representation of the finite group).  Two
of the motivating problems of the present survey are then the following:

\bigskip

{\bf Problems.}  i) Given a dimension $n$, determine the finite groups $G$ which admit
an action on a sphere or a homology sphere of dimension $n$ (comparing with the
class of finite subgroups of the orthogonal group ${\rm SO}(n+1)$).

\medskip

ii) Given a finite group $G$, determine the minimal dimension of a sphere or
homology sphere on which $G$ admits an action (show that it coincides with the
minimal dimension of a linear action of $G$ on a sphere).

\medskip

We close the introduction with a general result by Dotzel and Hamrick ([DH]),
again for finite $p$-groups: If $G$ is a finite $p$-group acting smoothly on a mod
$p$ homology $n$-sphere then $G$ admits also a linear action on $S^n$ such that the
two actions have the same dimension function for the fixed point sets of all
subgroups of $G$.

\medskip

In the next sections we will discuss finite groups acting on low-dimensional
spheres, starting with the 2-sphere $S^2$.

\bigskip
\vfill  \eject

{\bf 2.  Finite groups acting on $S^2$ and finite subgroups of ${\rm SO}(3)$.}

\medskip

Let $G$ be a finite group of (orientation-preserving) homeomorphisms of the
2-sphere $S^2$. It is a classical result of Brouwer and Kerekjarto from 1919
that such a finite group action on a 2-manifold is {\it locally linear}: each
fixed point of a nontrivial element has a regular neighbourhood which is a
2-disk on which the cyclic subgroup fixing the point acts as a standard
orthogonal rotation on the 2-disk. This implies easily that the quotient space
$S^2/G$ (the space of orbits) is a again a 2-manifold, or better a {\it
2-orbifold} of some signature $(g; n_1, \ldots, n_r)$:  an orientable  surface
of some genus $g$, with $r$ branch points of orders  $n_1, \ldots, n_r$ which
are the projections of the fixed points of  the nontrivial cyclic subgroup of
$G$ (of orders $n_i>1$).

\medskip

The projection $S^2  \to  S^2/G$ is a {\it branched covering}.  We choose some
triangulation of the quotient orbifold $S/G$ such that the branch points are
vertices of the triangulation, and lift this triangulation to a triangulation
of $S^2$; then the projection $S^2  \to  S^2/G$ becomes a simplicial map. If
the projection $p$ is a covering in the usual sense (unbranched, i.e. without
branch points), then clearly the Euler characteristic $\chi$ behaves
multiplicatively, i.e.  $2 = \chi (S^2) = |G| \; \chi(S^2/G) = |G| \; (2-2g)$
(just multiplying Euler characteristics with the order $|G|$ of $G$). In the
case with branch points, we can correct this by subtracting $|G|$ for each
branch point and then adding $|G|/n_i$ (the actual number of points of $S^2$
projecting to the i'th branch point), obtaining in this way the classical {\it
formula of Riemann-Hurwitz}:

$$2 = \chi (S^2) = |G| \; (2-2g - \sum_{i=1}^r (1 - {1\over n_i})).$$

It is easy to see that the only  solutions with positive integers of this equation  are
the following (first two columns):

\medskip

$(0; n,n), \hskip 10mm    |G| = n;  \hskip 10mm   G \cong \Bbb Z_n \hskip 3mm$ cyclic;

$(0; 2,2,n), \hskip 7mm   |G| = 2n; \hskip 8mm  G \cong \Bbb D_{2n} \hskip 2mm$ dihedral;

$(0; 2,3,3), \hskip 7mm   |G| = 12; \hskip 9mm   G \cong \Bbb \A_4  \hskip 3mm$
tetrahedral;

$(0; 2,3,4), \hskip 7mm   |G| = 24;  \hskip 9mm   G \cong \Bbb \S_4  \hskip 3mm$
octahedral;

$(0; 2,3,5), \hskip 7mm   |G| = 60;  \hskip 9mm  G \cong \Bbb \A_5 \hskip 3mm$
 dodecahedral.

\medskip

We still have to identify the groups $G$. For this, we first determine the
finite subgroups of the orthogonal group ${\rm SO}(3)$ and suppose that the
action of $G$ on $S^2$ is orthogonal. Of course the possibilities for the
signatures of the quotient orbifolds $S^2/G$ and the orders $|G|$ remain the
same, and for orthogonal actions one can identify the groups now as the
orientation-preserving symmetry groups of the platonic solids. As an example,
if the signature is  (0; 2,3,5) and  $|G| = 60$, one considers a fixed point
$P$ of a cyclic subgroup $\Z_5$ of $G$ and the five fixed points of subgroups
$\Z_3$ closest to $P$ on $S^2$; these are the vertices of a regular pentagon on
$S^2$ which is one of the twelve pentagons of a regular dodecahedron projected
to $S^2$, invariant under the action of $G$. Hence $G$, of order 60, coincides
with the orientation-preserving isometry group of the regular dodecahedron, the
dodecahedral group $\A_5$; see [W, section 2.6] for more details. In this way
one shows that the finite subgroups of ${\rm SO}(3)$ are, up to conjugation,
exactly the {\it polyhedral groups} as indicated in the list above: cyclic
$\Z_n$, dihedral $\D_{2n}$, tetrahedral $A_4$, octahedral $\S_4$ and
dodecahedral $A_5$.

\medskip

As a consequence, returning to topological actions of finite groups $G$ on $S^2$, the
topological orbifolds $S^2/G$ in the above list are geometric since they are
homeomorphic exactly to the quotients of $S^2$ by the polyhedral groups. Hence the
topological orbifolds $S^2/G$ have a spherical orbifold structure (a Riemannian
metric of constant curvature one, with singular cone points  of angles $2\pi/n_i$);
lifting this spherical structure to $S^2$ realizes
$S^2$ as a spherical manifold (i.e., with a Riemannian metric of constant curvature 1,
without singular points). Since the spherical metric on $S^2$ is unique up to isometry
this gives the standard Riemannian $S^2$, and
$G$ acts by isometries now. Hence every topological action of a finite group $G$ on
$S^2$ is conjugate to an orthogonal action, that is finite group actions on $S^2$ are
linear (this can be considered as the {\it orbifold geometrization in dimension two}, in
the spherical case).

\medskip

Concluding and summarizing, finite group actions on $S^2$ are {\it locally
linear}, and then also {\it linear}; the finite groups occurring are exactly
the polyhedral groups, and every topological action of such a group is
geometric (or linear, or orthogonal), i.e. conjugate to a linear action.

\bigskip

{\bf 3.  Finite groups acting on $S^3$ and finite subgroups of ${\rm SO}(4)$.}

\bigskip

{\bf 3.1.  Geometrization of finite group actions on $S^3$}

\medskip

We consider actions of a finite group $G$ on the 3-sphere $S^3$ now. The
first question is if such a topological action is {\it locally linear}; by the normal
form for orthogonal matrices, in  dimension three this means that an element with fixed
points acts locally as a standard rotation around some axis (the
orientation-preserving case). Unfortunately this is no longer true in dimension
three; after a first example of Bing from 1952 in the orientation-reversing case,
Montgomery-Zippen 1954 gave examples also of  orientation-preserving cyclic group
actions on $S^3$ with "wildly embedded fixed point sets", i.e. with fixed points
sets which are locally not homeomorphic to the standard embedding of $S^1$ in $S^3$.
Obviously such actions are not locally linear, and in particular cannot be conjugate
to smooth or linear actions.

\medskip

We will avoid these wild phenomena in the following by concentrating on {\it smooth} or
{\it locally linear} actions.  Suppose now that a cyclic group  $G \cong Z_p$, for a
prime $p$,  acts locally linear on $S^3$ with nonempty fixed point set; then it acts
locally as a standard rotation around an axis, and by compactness this axis closes
globally to an embedded knot $K \cong S^1$ in $S^3$. We note that also for a topological
action of $G$,  by general Smith fixed point theory the fixed point set of $G$ is a
knot $K$, i.e. an embedded  $S^1$ in $S^3$. If the action is locally linear, $K$
is a tame knot as considered in classical knot theory (smooth or polygonal
embeddings), otherwise $K$ is a wild knot leading to the different field of
wild or Bing topology (see [Ro, chapter 3.I] for such wild phenomena in
dimension three).

\medskip

So we will consider only {\it locally linear} actions in the following.  Globally, the
question arises then which knots $K$ can occur as the fixed point set of an action of a
cyclic group on $S^3$; it is easy to see that the action of  $G \cong Z_p$ is
{\it linear} (conjugate to a linear action) if and only if $K$ is a trivial knot
(unknotted, i.e.  bounding a disk in $S^3$). The classical {\it Smith conjecture}
states that $K$ is always a trivial knot, and that consequently locally linear actions
of cyclic groups are linear.  A positive  solution of the Smith conjecture was the first
major success of Thurston's geometrization program for 3-manifolds (see [MB]). This has
been widely generalized by Thurston then who showed that finite nonfree group actions on
closed 3-manifolds are build from geometric actions ({\it orbifold geometrization in
dimension three}), and recently by Perelman also for free actions of finite groups ({\it
manifold geometrization in dimension three}). As a consequence, every finite group acting
smoothly or locally linearly on $S^3$ is geometric, i.e. conjugate to an orthogonal or
linear action.

\medskip

Concluding, finite group actions on $S^3$ are not locally linear, in general, but
smooth or locally linear actions are linear; in particular, the finite groups
acting smoothly or locally linearly on $S^3$ are exactly the finite subgroups of
the orthogonal group ${\rm SO}(4)$.

\medskip

In the remaining part of this section we
discuss the finite subgroups of the orthogonal group ${\rm SO}(4)$, starting with the
relation between
${\rm SO}(3)$ and the unit quaternions.

\bigskip

{\bf 3.2.  The orthogonal group ${\rm SO}(3)$ and the unit quaternions $S^3$.}

\medskip

The orthogonal group ${\rm SO}(3)$ is homeomorphic to the real projective space
$\R\P^3$ of dimension three. In fact, by the normal form for orthogonal $3
\times 3$ matrices, such a matrix induces a clockwise rotation of the unit
3-ball in $\R^3$ around some oriented axis or diameter; parametrising the
diameter by an rotation angle from $-\pi$ to $\pi$, one obtains ${\rm SO}(3)$
by identifying diametral points on the boundary $S^2$ of the 3-ball (since
$-\pi$ and $\pi$  give the same rotation), and consequently  ${\rm SO}(3)$ is
homeomorphic to $\R\P^3$.

\medskip

Hence the universal covering of  ${\rm SO}(3)  \cong \R\P^3 \cong S^3/\langle
\pm {\rm id}\rangle$ is the 3-sphere $S^3$. Considering $S^3$ as the unit
quaternions, an orthogonal action of $S^3$ on the 2-sphere $S^2$  is obtained
as follows. The unit quaternions $S^3$ act on itself by conjugation $x \to
q^{-1}xq$, for a fixed $q \in S^3$; this action is clearly linear and also
orthogonal. Since $q$ fixes both poles 1 and -1 in $S^3$, it restricts to an
orthogonal action on the corresponding equatorial 2-sphere $S^2$ in $S^3$, so
this defines an element of the orthogonal group ${\rm SO}(3)$ and a group
homomorphism $S^3 \to {\rm SO}(3)$ of Lie groups of the same dimension, with
kernel $\pm 1$; by standard facts about Lie groups, this is the universal
covering of ${\rm SO}(3) \cong S^3/\langle \pm 1 \rangle$.

\medskip

The finite subgroups of ${\rm SO}(3)$ are exactly the {\it polyhedral groups}
$\Z_n,  \D_{2n}, \A_4,  \S_4$ and  $\A_5$. Their preimages in
the unit quaternions $S^3$ are the {\it binary polyhedral groups}
$$\Z_{2n}, \;  \D_{2n}^*, \; \A_4^*, \;  \S_4^*, \; \A_5^*$$
(cyclic, binary dihedral, binary tetrahedral, binary octahedral or binary
dodecahedral).  Since $S^3$ has a unique nontrivial involution  -1, together with the
cyclic groups of odd order these are exactly the finite subgroups of the unit
quaternions $S^3$, up to conjugation. By right multiplication, they act freely and
orthogonally on the 3-sphere
$S^3$, and the quotient spaces are examples of spherical 3-manifolds;
for example, $S^3/\A_5^*$ is the {\it Poincar\'e homology 3-sphere}.

\medskip

We note that the Lie group $S^3$ has various descriptions: it occurs as the unitary
group SU(2) over the complex numbers, as the universal covering group  Spin(3) of SO(3)
over the reals, and finally as the symplectic group Sp(1) over the quaternions (in
fact, the unit quaternions).

\bigskip

{\bf 3.3.  The orthogonal group ${\rm SO}(4)$ as a central product
$S^3 \times_{\Z_2}S^3$.}

\medskip

Passing to the orthogonal group ${\rm SO}(4)$ now acting on the unit 3-sphere
$S^3 \subset  \R^4$, there is an orthogonal action of $S^3 \times S^3$ on $S^3$
given by $x \to q_1^{-1}xq_2$, for a fixed pair of unit quaternions $(q_1,q_2) \in
S^3  \times S^3$. This defines again a homomorphism of Lie groups
$S^3 \times S^3 \to  {\rm SO}(4)$ of the same dimension, with kernel $\Z_2$ generated
by $(-1,-1)$, so this is in fact the universal covering of the Lie group ${\rm SO}(4)$.
In particular, the universal covering group ${\rm Spin}(4)$ of ${\rm SO}(4)$ is
isomorphic to
$S^3 \times S^3$, and ${\rm SO}(4)$ is isomorphic to the central product $S^3
\times_{\Z_2}S^3$ of two copies of the unit quaternions (the direct product with
identified centers, noting that the center of the unit quaternions is isomorphic
to $\Z_2$ generated by $-1$).

\medskip

Identifying ${\rm SO}(4)$ with $S^3 \times_{\Z_2}S^3$, the finite subgroups of
${\rm SO}(4)$ are, up to conjugation, exactly the finite subgroups of the central
products
$$P_1^* \times_{\Z_2} P_2^* \subset S^3 \times_{\Z_2} S^3$$  of two binary polyhedral
groups $P_1^*$ and $P_2^*$. The most interesting example of such a group is the
central product $\A_5^* \times_{\Z_2}\A_5^*$ of two binary dodecahedral groups
which occurs as the orientation-preserving symmetry group of the regular
4-dimensional 120-cell (a fundamental domain for the universal covering group
$\A_5^*$ of the Poincar\'e homology sphere $S^3/\A_5^*$ is a regular spherical
dodecahedron, and 120 copies of this dodecahedron give a regular spherical
tesselation of the 3-sphere $S^3$;  the vertices of this tesselation are the
vertices of the regular euclidean 120-cell in $\R^4$ whose faces are 120
regular dodecahedra.)

\medskip

Concluding, the finite subgroups of ${\rm SO}(4)$ are exactly the subgroups of the
central products $P_1^* \times P_2^*$  of two binary polyhedral groups $P_1^*$ and
$P_2^*$. It is then an algebraic exercise to classify the possible groups, up to
isomorphism and up to conjugation (see [DV] for a list of the finite subgroups of
${\rm SO}(4)$, and also of ${\rm O}(4)$).

\bigskip

{\bf 4.  Finite groups acting on $S^4$ and finite subgroups of ${\rm SO}(5)$.}

\medskip

As we have seen in the previous sections, finite group actions on the 2-sphere
are locally linear, and then also linear. In dimension three, finite group
actions are not locally linear, in general, but by deep results of Thurston and
Perelman, smooth or locally linear actions on the 3-sphere are linear. In
dimension four, also smooth or locally linear actions are no longer linear, in
general. In fact it has been shown by Giffen in 1966 that the Smith conjecture
fails in dimension four, by constructing examples of smooth actions of a finite
cyclic group on $S^4$ whose fixed point sets are  knotted 2-spheres in  $S^4$
(see [R, chapter 11.C]); in particular, such an action cannot be linear.

\medskip

Restricting again to smooth or locally linear actions, we consider now the
problem stated in the introduction: which finite groups $G$ admit a smooth
orientation-preserving action  on the 4-sphere $S^4$; also, what are the finite
subgroups of ${\rm SO}(5)$?

\medskip

Suppose that $G$ is a finite group with an orientation-preserving, faithful, linear
action on $S^4 \subset \R^5$; using the language of group representations, this
means that $G$ has a faithful, orientation-preserving, real representation in
dimension five.  If such a representation is {\it reducible} (a direct sum of
lower-dimensional representations), $G$ is an orientation-preserving subgroup of a
product of orthogonal groups ${\rm O}(3) \times {\rm O}(2)$ or
${\rm O}(4) \times {\rm O}(1)$, so one can reduce to
lower dimensions.

\medskip

Suppose that the representation is {\it irreducible} but {\it imprimitive};
this means that there is a decomposition of $\R^5$ into proper linear subspaces which
are permuted transitively by the group. Since the dimension five is prime, these
linear subspaces have to be 1-dimensional (such a representation is then called
monomial). The group of orthogonal maps permuting the five factors $\R$ of $\R^5$ is
the Weyl-group  $W_5 = (\Z_2)^5 \rtimes \S_5$ of inversions and permutations of
coordinates, i.e. the semidirect product of the normal subgroup $(\Z_2)^5$
generated by the inversions and the symmetric group $\S_5$ of permutations of the
factors.  Hence $G$ is  a subgroup of the Weyl-group $W = (\Z_2)^5 \rtimes \S_5$ in this
case.

\medskip

There remains the case of an {\it irreducible, primitive} representation; this
is the main case which has been considered by various authors and for arbitrary
dimension; a major problem here is to find the simple groups which admit such a
representation (i.e.,  groups without a nontrivial proper normal subgroup).
This leads into classical representation theory of finite groups, and we will
not go further into it. In fact, we gave the above description mainly as a
motivation for the next result on smooth or locally linear actions of finite
groups on the 4-sphere.

\bigskip
\vfill  \eject

{\bf  Theorem 1.} ([MeZ1])  {\sl  A finite group $G$ with a smooth or locally
linear, orientation-preserving action on the 4-sphere, or on any homology 4-sphere, is
isomorphic to one of the following groups:

\smallskip

i)  an orientation-preserving subgroup of ${\rm O}(3) \times {\rm O}(2)$ or
${\rm O}(4) \times {\rm O}(1)$;

\smallskip

ii) an orientation-preserving subgroup of the Weyl group $W = (\Z_2)^5 \rtimes \S_5$;

\smallskip

iii)  $\A_5$, $\S_5$, $\A_6$ or $\S_6$;

\smallskip

i') if $G$ is nonsolvable, a 2-fold extension of a subgroup of ${\rm SO}(4)$.}

\medskip

Note that the different cases of the Theorem are not mutually exclusive. The
only indetermination remains case i'); in fact, in this case $G$ should be
isomorphic to a subgroup of ${\rm O}(4)$ and hence to an orientation-preserving
subgroup  of the group ${\rm O}(4) \times {\rm O}(1)$ of case i); however the
proof in this case is not completed at present since many different cases have
to be considered, according to the long list of  finite subgroups of ${\rm
SO}(4)$ (see [DV]).

\bigskip

{\bf  Corollary 1.} {\sl  A finite group $G$ which admits an orientation-preserving
action on a homology 4-sphere is isomorphic to a subgroup of ${\rm SO}(5)$ or, if
$G$ is solvable, to a 2-fold extension of a subgroup of ${\rm SO}(4)$. }

\medskip

The symmetric group $\S_6$ acts orthogonally on $\R^6$ by permutation of
coordinates, and also on its subspace $\R^5$ defined by setting the sum of the
coordinates equal to zero (this is called the standard representation of
$\S_6$), and hence on the unit sphere $S^4 \subset \R^5$. Composing the
orientation-reversing elements by -id, one obtains an orientation-presering
action of $\S_6$ on $S^4$ (alternatively, $\S_ 6$ acts on the 5-simplex by
permuting its six vertices, and hence on its boundary which is the 4-sphere.)

\medskip

For linear action, Theorem 1 and its proof easily give the following
characterization of the finite subgroups of ${\rm SO}(5)$.

\bigskip

{\bf Corollary 2.} {\sl  Let $G$ be a finite subgroup of the orthogonal group
${\rm SO}(5)$.  Then one of the following cases occurs:

\smallskip

i) $G$ is conjugate to an orientation-preserving subgroup of  ${\rm O}(4)\times{\rm
O}(1)$ or ${\rm O}(3)\times{\rm O}(2)$ (the reducible case);

\smallskip

ii) $G$ is conjugate to  a subgroup of the Weyl group $W = (\Z_2)^4 \rtimes \S_5$
(the irreducible, imprimitive case);

\smallskip

iii) $G$ is isomorphic to $\A_5$, $\S_5$, $\A_6$ or $\S_6$  (the irreducible,
primitive case).}

\medskip

See the character tables in [C] or [FH] for the irreducible representations of
the groups in iii) (e.g. $\A_5$ occurs as an irreducible subgroup of all three
orthogonal groups ${\rm SO}(3)$, ${\rm SO}(4)$ and ${\rm SO}(5)$).

\medskip

It should be noted that the proof of Corollary 2 is considerably easier than
the proof of Theorem 1. For both Theorem 1 and Corollary 2 one has to determine
the finite simple groups which act on a homology 4-sphere resp. which admit an
orthogonal action on the 4-sphere. In the case of Theorem 1 this is based on
[MeZ2, Theorem 1] which employs the Gorenstein-Harada classification of the
finite simple groups of sectional 2-rank at most four (see [Su2], [G1]). For the proof of
Corollary 2 instead, this heavy machinery from the classification of the finite simple
groups can be replaced by much shorter arguments from the representation theory of
finite groups.

\medskip

For solvable groups $G$ instead, the proof of Theorem 1 is easier; here one
can consider the Fitting subgroup of $G$, the maximal normal nilpotent subgroup,
which is nontrivial for solvable groups. As a nilpotent group, the  Fitting
subgroup is the direct product of its Sylow $p$-subgroups, has nontrivial center and
hence nontrivial cyclic normal subgroups of prime order. A starting point of the
proof of Theorem 1 is then the following lemma which shows some of the basic ideas
involved.

\bigskip

{\bf Lemma 1.}  {\sl Let $G$ be a finite group with a
smooth, orientation-preserving action on a homology 4-sphere. Suppose that  $G$
has  a cyclic normal group $\Z_p$ of prime order $p$; by Smith fixed point
theory, the fixed point set of $\Z_p$ is either a 0-sphere $S^0$ or
a 2-sphere $S^2$ (i.e., a mod $p$ homology sphere of even codimension).

\smallskip

i) If the fixed point set of $\Z_p$ is a 0-sphere then $G$ contains of index at
most two a subgroup isomorphic to a subgroup of ${\rm SO}(4)$. Moreover if $G$
acts orthogonally on $S^4$ then $G$ is conjugate  to a subgroup of ${\rm
O}(4)\times {\rm O}(1)$.

\smallskip

ii) If the fixed point set of $\Z_p$ is a 2-sphere then $G$ is isomorphic to a
subgroup of ${\rm O}(3) \times {\rm O}(2)$. Moreover if $G$ acts orthogonally
on $S^4$ then $G$ is conjugate to a subgroup of ${\rm O}(3) \times {\rm
O}(2)$.}

\medskip

{\it Proof.} i)  Since $\Z_p$ is normal in $G$, the group $G$ leaves invariant the
fixed point set $S^0$ of $\Z_p$ which consists of two points. A subgroup $G_0$ of
index at most two of $G$ fixes both points and acts orthogonally and
orientation-preservingly on a 3-sphere, the boundary of a  $G_0$-invariant
regular neighborhood of one of the two fixed points.

If the action of $G$ is an
orthogonal action on the 4-sphere then $G$ acts orthogonally on the equatorial
3-sphere of the 0-sphere $S^0$ and hence is a subgroup of ${\rm O}(4)\times {\rm
O}(1)$, up to conjugation.

\medskip

ii) The group $G$ leaves invariant the fixed point set $S^2$ of $\Z_p$. A
$G$-invariant regular neighbourhood of $S^2$ is diffeomorphic to the product of
$S^2$ with a 2-disk, so $G$ acts on its boundary $S^2 \times S^1$ (preserving
its fibration by circles). Now, by the geometrization of finite group actions
in dimension three, it is well-known that every finite group action on $S^2
\times S^1$ preserves the product structure and is standard, i.e.  is conjugate
to a subgroup of its isometry group ${\rm O}(3) \times {\rm O}(2)$.

If $G$ acts orthogonally on $S^4$ then the group
$G$ leaves invariant $S^2$, the corresponding 3-dimensional subspace in $\Bbb{R}^5$
as well as its orthogonal complement, so up to conjugation it is a subgroup of ${\rm
O}(3) \times {\rm O}(2)$.

\bigskip

{\bf 5.  Higher dimensions}

\medskip

Relevant in the context of linear actions on spheres is the classical {\it Jordan
number}: for each dimension $n$ there is an integer $j(n)$ such that each finite
subgroup of the complex linear group ${\rm GL}_n(\Bbb C)$, and hence in particular
also of its subgroup ${\rm SO}(n)$, has a normal abelian subgroup of index at most
$j(n)$ (we note that a lower bound for $j(n)$ is $(n+1)!$ since the symmetric group
$\S_{n+1}$ is a subgroup of
${\rm GL}_n(\Bbb C)$; see the comments to [KP, Theorem 5.1] for an upper bound.)
Whereas this is insignificant for abelian group, it implies that the order of a
nonabelian simple groups acting linearly on
$S^n$ is bounded by
$j(n+1)$; in particular, up to isomorphism there are only finitely many finite
simple groups (always understood to be nonabelian in the following) which admit a
faithful, linear action on $S^n$ (or equivalently, have a faithful, real, linear
representation in dimension $n+1$).  For smooth or locally linear actions of finite
simple groups on spheres and homology spheres, there is the following analogue.

\bigskip

{\bf Theorem 2.} ([GZ]) {\sl For each dimension $n$, up to isomorphism there are
only finitely many finite simple groups  which admit
a smooth or locally linear action on the $n$-sphere, or on some homology sphere of
dimension  $n$.}

\medskip

We note that any finite simple group admits many smooth actions on high-dimensional
spheres which are not linear  (conjugate to a linear action; see the survey [Da, section
7]).

\medskip

It is natural to ask whether the Jordan number theorem can be  generalized for
all finite groups acting on homology $n$-spheres.  Since, as noted in the
introduction, finite  $p$-groups admitting a smooth action on some homology
$n$-sphere admit also a linear action on $S^n$ ([DH]), it is easy to generalize the
Jordan number theorem for nilpotent groups; so the Jordan theorem remains true for the
two extreme opposite cases of nilpotent groups and simple groups,  but at present we
don't know it for arbitrary finite groups.

\medskip

Not surprisingly, the proof of Theorem 2 requires the full classification of
the finite simple groups; we will present part of the proof in the following.
We note that the proof of Theorem 2 permits to produce for each dimension $n$ a
finite list of finite simple groups which are the candidates for actions on
homology $n$-spheres; then one can identify those groups from the list which
admit a linear action on $S^n$ (or equivalently, have a faithful, real, linear
representation in dimension $n+1$), and try to eliminate the remaining ones by
refined methods. For example, it is shown in [MeZ2-4] that the only finite
simple group which admits an action on a homology 3-sphere is the alternating
group $\A_5$, and that the only finite simple groups acting on a homology
4-sphere are the alternating groups $\A_5$ and $\A_6$; already these
low-dimensional results require heavy machinery from the classification of the
finite simple groups.

\medskip

Crucial for the proofs of Theorems 1 and 2 is a control over the  minimal
dimension of an action of a linear fractional group ${\rm PSL}_2(p)$ and a
linear group ${\rm SL}_2(p)$ (the latter is the group of $2Ê\times 2$-matrices
of determinant one over the finite field with $p$ elements, the former its
factor group by the central subgroup $\Z_2 = \langle \pm E_2 \rangle$);  this
is given by the following:

\bigskip

{\bf Proposition 1.} ([GZ]) {\sl For a prime $p \ge 5$, the minimal dimension
of an action of a linear fractional group ${\rm PSL}_2(p)$  on a mod $p$
homology sphere  is $(p-1)/2$ if  $p  \equiv  1  \; {\rm mod} \; 4$, and $p-2$
if $p \equiv 3 \; {\rm mod} \;  4$, and these are also lower bounds for the
dimension of such an action of a linear group  ${\rm SL}_2(p)$.}

\medskip

Whereas the groups ${\rm PSL}_2(p)$ admit linear actions on spheres of the
corresponding dimensions (see e.g. [FH]), for the groups ${\rm SL}_2(p)$ the
minimal dimension of a linear action on a sphere is $p-2$ resp. $p$,  that is
strictly larger than the lower bounds given in Proposition 1, so the minimal
dimension of an action on a homology sphere remains open here.

\medskip

In analogy with Proposition 1, we need also the following result from Smith
fixed point theory for  elementary abelian $p$-groups ([Sm]).

\bigskip

{\bf Proposition 2.}  {\sl  The minimal dimension of a faithful,
orientation-preserving action of an elementary abelian $p$-group $(\Bbb Z_p)^k$
on a mod $p$ homology sphere is $k$ if $p = 2$, and $2k-1$ if $p$ is an odd
prime.}

\medskip

Considering  commutator subgroups, the proof of the following lemma is an easy
exercise (a finite central extension of $G$ is a finite group with a
central subgroup whose factor group is isomorphic to $G$; a group is perfect if it
coincides with its commutator subgroup or, equivalently, the abelianized group is
trivial).

\bigskip

{\bf Lemma 2.} {\sl If a finite group $G$ has a perfect subgroup $H$ then any finite
central extension of $G$ contains a perfect central extension of $H$.}

\medskip

For the proof of Theorem 2, Lemma 2 will be applied mainly when $H$ is a linear
group ${\rm SL}_2(q)$, for a prime power $q = p^k$ (that is, over the finite field
with $p^k$ elements); we note that, for $q \ge 5$ and different from 9, the only
perfect central extension of the perfect group ${\rm SL}_2(q)$ is the group itself
(see [H, chapter V.25]) (and the only nontrivial perfect central extension of ${\rm
PSL}_2(q)$ is ${\rm SL}_2(q)$).

\bigskip

On the basis of Propositions 1 and 2 and the classification of the finite
simple groups, we indicate now the

\bigskip

{\it Proof of Theorem 2.}  Fixing a dimension $n$, we have to exclude all but
finitely many finite simple groups. By the classification of the finite simple
groups, a finite simple group is one of 26 sporadic groups, or an alternating
group, or a group of Lie type ([Co], [G1]). We can neglect the sporadic groups
and have to exclude all but finitely many groups of the infinite series.
Clearly, an alternating group $\A_m$ contains elementary abelian subgroups
$(\Bbb Z_2)^k$ of rank $k$ growing with the degree $m$, so Proposition 2
excludes all but finitely many alternating groups and we are left with the
infinite series of groups of Lie type.

\medskip

We consider first the projective linear groups ${\rm PSL}_m(q)$, for a prime
power $q=p^k$. The group ${\rm PSL}_2(q)$ has subgroups ${\rm PSL}_2(p)$ and
$(\Z_p)^k$ (the subgroup represented by all diagonal matrices with entries one
on the diagonal, isomorphic to the additive group of the field with $p^k$
elements), and by Propositions 1 and 2 only finitely many primes $p$ and prime
powers $p^k$ can occur. If $m \ge 3$ instead, ${\rm PSL}_m(q)$ has subgroups
${\rm SL}_2(q)$ and ${\rm SL}_2(p)$; again Proposition 1 excludes all but
finitely many primes $p$ and, since also ${\rm SL}_2(q)$ has an  elementary
abelian subgroup $(\Bbb Z_p)^k$, by Proposition 4 only finitely many powers
$p^k$ of a fixed prime $p$ can occur. Concluding, only finitely many prime
powers $p^k$ can occur for a fixed dimension $n$.  Note that, in a similar way,
also for the groups ${\rm SL}_2(q)$ only finitely many values of $q$ can occur.
We still have to bound $m$; if $q$ is not a power of two, then ${\rm PSL}_m(q)$
has a subgroup $(\Z_2)^{m-2}$ represented by diagonal matrices with entries
$\pm 1$ on the diagonal, so $m$ is bounded by Proposition 2. If $q$ is a power
of two, one may consider instead subgroups $(Z_p)^{[m/2]}  <  {\rm
SL}_2(p)^{[m/2]} < {\rm PSL}_m(q)$ and again apply Proposition 2.

\medskip

This finishes the proof of Theorem 2 for the case of the projective linear groups
${\rm PSL}_m(q)$.

\medskip

The proof for the unitary groups ${\rm PSU}_m(q)$ and the symplectic groups  ${\rm
PSp}_{2m}(q)$ is similar, noting that there are isomorphisms
$${\rm PSU}_2(q) \cong {\rm PSp}_2(q) \cong {\rm PSL}_2(q), \;\; {\rm SU}_2(q) \cong
{\rm Sp}_2(q) \cong {\rm SL}_2(q);$$ in particular, if $m \ge 3$ resp. $2m
\ge 4$, the latter groups are subgroups of both ${\rm PSU}_m(q)$ and ${\rm
PSp}_{2m}(q)$, so we can conclude as before.

\medskip

The last class of classical groups are the orthogonal groups
$\Omega_{2m+1}(q) = P\Omega_{2m+1}(q)$ and $P\Omega^\pm_{2m}(q)$ (the latter stands
for two different groups which are simple if $m \ge 3$). There are isomorphisms
$$\Omega_3(q) \cong {\rm PSL}_2(q), \; P\Omega^+_4(q) \cong {\rm PSL}_2(q) \times
{\rm PSL}_2(q), \; P\Omega^-_4(q) \cong {\rm PSL}_2(q^2)$$ (see [Su, p.384]).
By canonical inclusions between orthogonal groups and the cases considered before,
this leaves again only finitely many possibilities.

\medskip

Next we consider the exceptional groups $G_2(q)$, $F_4(q)$, $E_6(q)$, $E_7(q)$,
$E_8(q)$ as well as the Steinberg triality groups  $^3D_4(q)$.  By [St, Table 0A8],
[GL, Table 4-1], up to central extensions there are inclusions
$$E_6(q)  >  F_4(q) >  {^3D_4(q)}  > G_2(q) > {\rm PSL}_3(q),
\;\; E_7(q) > {\rm PSL}_8(q), \;\;  E_8(q) > {\rm PSL}_9(q).$$ Applying Lemma 1 we
reduce to subgroups ${\rm SL}_2(q)$ in all cases.

\bigskip

Finally, there remain the twisted groups
$$^2E_6(q), \;\;  ^3D_4(q), \;\;  Sz(2^{2m+1}) = {^2B_2(2^{2m+1})},\;
\; ^2G_2(3^{2m+1}), \;\;  ^2F_4(2^{2m+1}).$$ By [St, Table 0A8], [GL, Table 4-1],
up to central extensions there are  inclusions $^2E_6(q) > F_4(q)$ and  $^3D_4(q) >
G_2(q)$ (already considered before); the Suzuki groups $Sz(2^{2m+1})$ have subgroups
$(\Z_2)^{2m+1}$ ([G1, p.74]), the Ree groups $^2G_2(3^{2m+1})$ subgroups
${\rm PSL}_2(3^{2m+1})$ ([G2, p.164]) and the Ree groups
$^2F_4(2^{2m+1})$ subgroups ${\rm SU}_3(2^{2m+1})$ ([GL, Table 4-1]), so in all
these cases some  previously considered case applies.

\medskip

This completes the proof of Theorem 2.

\bigskip \bigskip

{\bf Acknowledgment.}  I want to thank the colleagues  and friends from the
Universities of Novosibirsk, Kemerovo and Krasnojarsk for their great
hospitality during my visit in june 2011 which gave origin to the present
notes.

\bigskip \bigskip

\centerline {\bf References}

\bigskip

\item {[Bre]} G. Bredon, {\it Introduction to Compact Transformation Groups.}
Academic Press, New York 1972

\item {[Bro]} K.S. Brown, {\it Cohomology of Groups.}  Graduate Texts in Mathematics
87, Springer 1982

\item {[C]} J.H. Conway, R.T. Curtis, S.P. Norton, R.A. Parker, R.A. Wilson, {\it
Atlas of Finite Groups.} Oxford University Press 1985

\item {[Da]}  M.W. Davis, {\it  A survey of results in higher dimensions.}  The Smith
Conjecture, edited by J.W. Morgan, H. Bass, Academic Press 1984,  227-240

\item {[DH]}  R.M.Dotzel, G.C.Hamrick, {\it  $p$-group actions on homology
spheres.}  Invent. math. 62, 437-442  (1981)

\item {[DV]}  P. Du Val, {\it  Homographies, Quaternions and Rotations.} Oxford
Math. Monographs, Oxford University Press 1964

\item {[FH]} W. Fulton, J. Harris, {\it Representation Theory: A First Course.}
Graduate Texts in Mathematics 129,  Springer-Verlag 1991

\item {[G1]} D. Gorenstein, {\it The Classification of Finite Simple Groups.} Plenum
Press, New York 1983

\item {[G2]} D. Gorenstein, {\it Finite Simple Groups: An
Introduction to  their Classification.}  Plenum Press, New York
1982

\item {[GL]} D. Gorenstein, R. Lyons, {\it The local structure of
finite group of characteristic 2 type.} Memoirs Amer. Math. Soc.
276, vol. 42, 1-731 (1983)

\item {[GZ]} A. Guazzi, B. Zimmermann, {\it On finite simple groups acting on
homology spheres.}   Pre-print Trieste 2011 (see arXiv:1106.1067  for  a preliminary
version)

\item {[H]} B. Huppert, {\it Endliche Gruppen I.}  Grundlehren der
math. Wissenschaften, Band 134, Berlin: Springer-Verlag 1967

\item {[KP]} J. Kuzmanovich, A. Pavlichenko,  {\it Finite groups of matrices whose
entries are integers.}  Monthly 109,  173-186  (February 2002)

\item {[MeZ1]} M. Mecchia, B. Zimmermann, {\it On finite  groups acting on homology
4-spheres and finite subgroups of ${\rm SO}(5)$.}  Top. Appl. 158,  741-747  (2011)

\item {[MeZ2]} M. Mecchia, B. Zimmermann, {\it On finite simple and nonsolvable
groups acting on homology 4-spheres.}  Top. Appl. 153, 2933-2942 (2006)

\item {[MeZ3]} M. Mecchia, B. Zimmermann, {\it On finite simple groups acting on
integer and mod 2 homology 3-spheres.}  J. Algebra 298, 460-467  (2006)

\item {[MeZ4]} M. Mecchia, B. Zimmermann, {\it On finite simple groups acting on
integer and mod 2 homology 3-spheres.}  J. Algebra 298, 460-467  (2006)

\item {[MB]} J.Morgan, H.Bass, {\it The Smith Conjecture.} Academic Press, New York
1984

\item {[R]} D. Rolfson, {\it Knots and Links.}  Mathematics Lecture Series 7,
Publish of Perish,  Berkeley 1976

\item {[St]} E. Stensholt, {\it Certain embeddings among finite
groups of Lie type.}  J. Algebra 53, 136-187  (1978)

\item {[Su1]} M. Suzuki, {\it Group Theory II.}  Springer-Verlag 1982

\item {[Su2]} M. Suzuki, {\it Group Theory II.}  Springer-Verlag 1982

\item {[W]} J.A. Wolf, {\it Spaces of Constant Curvature.}

\bye